\documentclass{tran-l}

\numberwithin{equation}{section}
\theoremstyle{plain}

 \newtheorem{thm}{Theorem}[]

 \newtheorem{lem}[thm]{Lemma}
 \newtheorem{prop}[thm]{Proposition}

 \theoremstyle{definition}
 \newtheorem{propers}[]{Properties}
 \newtheorem{defn}[]{Definition}
 \newtheorem{defns}{Definitions}[]
 \newtheorem{rem}{Remark}[]

\begin{document}

\title[]
 {An explicit Construction of the Jacobian of a
nonsingular curve and its group law}

\author{ Tal Perri }

%\date{\today}

\address{Department of Mathematics \& Statistics, Bar-Ilan University, Ramat-Gan, Israel}

\email{talperri@zahav.net.il}

%%% ----------------------------------------------------------------------

\begin{abstract}
This paper is devoted to constructing an explicit efficient
representation for the Jacobian variety of a
nonsingular curve of genus greater than 1, and its group law.\\
We start by introducing some theoretic background for the theme in
issue, based on a paper of G.W.Anderson.\\
Following this theory we describe an algorithm for executing the
group law on the Jacobian elements, consisting of several steps.\\
Finally we exhibit two examples. The first example exhibits an
explicit formulation for the Jacobian of a well known non singular
elliptic curve of genus $1$, that is, the Weierstrass curve. The
second example introduces an explicit computation of the group law
for a non singular curve of genus 3,
$$C :  \ \omega_{2}^{4}=\omega_{1}(\omega_{1}-\omega_{0}))(\omega_{1}-2\omega_{0})(\omega_{1}-3\omega_{0})$$
where $\omega_{0}, \omega_{1}, \omega_{2}$ are homogeneous
coordinates of the projective plane  $\mathbb{P}^{ \ 2}$. \\
These examples are highly detailed for the convenience of the
reader.\\
\end{abstract}
%%% ----------------------------------------------------------------------
\maketitle
%%% ----------------------------------------------------------------------
\section*{introduction}
This paper is devoted to constructing an explicit efficient
representation for the Jacobian variety of a nonsingular curve of
genus greater than 1, and its group law.\\
The construction of the representation for the Jacobian variety
exhibited in sections one to three is mainly taken from Anderson's
paper \cite{And02}, where in Lemma 10 we decompose the process into
well organized steps to be computed, thus giving it an algorithmic
structure . The example introduced in the fifth section follows
these exact steps to reach the desired result.\\
The first example given in section four appears partly in Anderson's
paper. Our approach here is different and much more detailed. The
second example (and
the most interesting one) is new.\\
 The construction includes three stages. The first is a
construction of a matrix representation for the equivalence classes
of divisors, where each  divisor class defines a special class of
matrices, referred to as a $G-form$, representing it uniquely. There
is a bijective correspondence between the divisor classes and the
$G-form$ classes. The second construction is that of the Jacobian's
group law implemented on these $G-forms$ classes. The third
construction is that of the Abel map which maps $G-form$ classes
bijectively to a set of proportionality classes of special kind of
matrices named $Jacobi$ matrices. The set of $Jacobi$ matrices forms
a projective algebraic variety.\\

By constraints of the extent of this paper, we decide to leave the
construction of the inverse Abel map out of the scope of this paper.

\newpage
\section*{preliminaries}
\text{}\\
Throughout this paper, unless mentioned otherwise, we take :\\
$\bullet$ $k$ an algebraically closed field. \\
$\bullet$ $n$ an integer $\geq 2$.\\
$\bullet$ $A$ a finitely generated $k$-algebra without zero divisors.\\
$\bullet$ $L \subset A$ a finite dimensional $k$-subspace. \\
$\bullet$ $k$-\emph{algebras} are taken to be commutative with unit,
where elements of $k$ are often \text{ \ } referred to as
\emph{scalars} or \emph{constants}.

\section*{the abeliant and the discriminant}
\begin{defn}
Let $\{X^{(l)}\}_{l=0}^{n+1}$ be a family of $n$ by $n$ matrices
with entries in a ring $R$ and let $\{s_{i}\}_{i=1}^{n} \cup
\{t_{j}\}_{j=1}^{n}$ be a family of independent variables.\\

The $\mathbf{abeliant}$ ,denoted by $abel(X^{0},...,X^{n+1})$ or
$abel_{l=0}^{n+1}X^{(l)}$, is defined to be an $n$ by $n$ matrix the
entry of which in position $ij $ is the coefficient with which the
monomial\\

\text{ \ } \ \text{ \ } \ \text{ \ } \ $ s_{i}^{-1}t_{j}^{-1} \cdot
\prod_{a=1}^{n}{s_{a}} \cdot \prod_{b=1}^{n}{t_{b}} = s_{1}\cdots
\hat{s_{i}} \cdots
s_{n}t_{1}\cdots \hat{t_{j}} \cdots t_{n} $\\

appears in the expansion of the expression\\
\text{ \ } \ \text{ \ } \ $
trace(X^{(0)}(\sum_{b=1}^{n}{t_{b}X^{(b)}})^{\star}X^{(n+1)}(\sum_{a=1}^{n}{s_{a}X^{(a)}})^{\star})
$\\
\text{ \ } \\
as an $R$-linear combination of monomials in the $s$'s and the
$t$'s.\\

The $\mathbf{discriminant}$ $\bigtriangleup(X^{1},...,X^{n+1})$
(also denoted by $\bigtriangleup_{l=1}^{n+1}X^{(l)}$), of the given
family of matrices is defined to be
$$ | \sum_{l=1}^{n}{X^{(l)}} |^{ \ 2n-2} \cdot \prod_{i=1}^{n}{|\sum_{l \in \{1,...,n+1 \}
\setminus \{i\}}{X^{(l)}}|^{\ 2}}$$
\end{defn}
%\begin{propers}{\textbf{}} \\
 %Let $\{X^{(l)}\}_{l=0}^{n+1}$ be a family of $n$ by $n$
%matrices
%with entries in a ring $R$. \\
%Let $U , V$ be n by n matrices with entries in $R$.\\
%Let $\langle 0 | n+1 \rangle$ denote the permutation of $
%\{0,...,n+1 \}$ exchanging $0$ and $n+1$ and fixing all other
%elements ,and let $X^{T}$ denote the transpose of a square matrix
%$X$.\\

%The following equalities hold:\\

%$\bullet$ \ $abel_{l=0}^{n+1}(UX^{(l)}V)=|U|^{2}|V|^{2}abel_{l=0}^{n+1}X^{(l)}$\\

%$\bullet$ \ $abel_{l=0}^{n+1}X^{(\langle 0 | n+1 \rangle
%l)}=(abel_{l=0}^{n+1}X^{(l)})^{T} = abel_{l=0}^{n+1}(X^{(l)})^T $\\

%$\bullet$ \ $\bigtriangleup_{l=1}^{n+1}{UX^{(l)}V}=|U|^{ \
%4n-2}|V|^{
%\ 4n-2}\bigtriangleup_{l=1}^{n+1}{X^{(l)}}$\\

%$\bullet$ \
%$\bigtriangleup_{l=1}^{n+1}{X^{(l)}}=\bigtriangleup_{l=1}^{n+1}{(X^{(l)})^{T}}$\\

Various properties of these operators can be found in
Anderson's paper \cite{And02}.\\

%\end{propers}

%\begin{cor}{\textbf{The key relations}} \\
%Let $\{X^{(l)}\}_{l=0}^{n+1}$ be a family of $n$ by $n$ matrices
%such that for each $l=1,...n+1$ there exists a factorization
%$$X^{(l)}=u^{(l)}v^{(l)}$$
%where $u^{(l)}$ and $v^{(l)}$ are column and  row vectors
%respectively, with entries in $R$. Then the following equalities
%hold
%\begin{equation}\label{Key1}
%abel_{l=0}^{n+1}X^{(l)}=MU^{\star}X^{(0)}V^{\star}L
%\end{equation}

%\begin{equation}\label{Key2}
%\bigtriangleup_{l=1}^{n+1}X^{(l)}=|MU^{\star}|^{\ 2}|V^{\star}L|^{\
%2}
%\end{equation}

%where
 %$$U=\begin{pmatrix}  u^{(1)} & \cdots & u^{(n)}
%\\ \end{pmatrix}, \ \ V=\begin{pmatrix} & v^{(1)}\\ & \vdots \\ & v^{(n)}
%\\ \end{pmatrix}$$
%$$M=diag\begin{pmatrix}  (v^{(n+1)}V^{\star})_{1},  \cdots,
%(v^{(n+1)}V^{\star})_{n} \\\end{pmatrix}, \  L=diag\begin{pmatrix}
%(U^{\star}u^{(n+1)})_{1}, \cdots,
% (U^{\star}u^{(n+1)})_{n} \\\end{pmatrix}$$
%\end{cor}
%\text{ \ } \\

\section*{segre matrices and jacobi matrices}

Let $X$ and $Y$ be matrices with entries in a $k$-algebra $R$.

\begin{defns}{\textbf{}} \\
$\bullet$ $X$ is said to be \emph{k-general} if there exist a row of
$X$ and a column of $X$ each \text{ \ } of which consists of
k-linearly
independent entries.\\
$\bullet$ $X$ and $Y$ are said to be \emph{k-equivalent} if there
exist square matrices $\Phi$ and $\Psi$ with \text{ \ } entries in
$k$ such that $|\Phi| \neq 0 \ , \ |\Psi| \neq 0$ such that $Y=\Phi
X
\Psi$.\\
$\bullet$ Let $x$ and $y$ be vectors in a common vector space over
$k$, $x$ is said to be \\
\text{ \ } \emph{k-proportional} to $y$ if there exists a nonzero
scalar $c \in k$ such that $x=cy$.
\end{defns}
\goodbreak
\begin{defn}
Recalling the notations introduced in the Preliminaries, a
\textbf{Segre matrix} $X$ is an $n$ by $n$ matrix with entries in
$L$, such that $X$ is k-general and $rank(X) \leq 1$ .
\end{defn}
We say $X$ is of type $(k,n,A,L)$
\goodbreak
\begin{propers}{\textbf{}} \\
$\bullet$ Let $X$ be a Segre matrix, then any matrix with entries in
$A$ to which $X$ is\\
 \text{ \ } k-equivalent is a Segre matrix.\\
$\bullet$ The transpose $X^{T}$ is a Segre matrix.\\
$\bullet$ There exists a factorization $X=uv$ where $u$ (resp. $v$)
is a column (resp. row)\\
 \text{ \ } vector with entries in the
fraction field of
$A$.\\
$\bullet$ Given any such factorization $X=uv$ ,the entries of $u$
(resp. $v$)are k-linearly\\
 \text{ \ } independent.\\
$\bullet$ Given any two such factorizations $X = uv =u'v'$, there
exists unique nonzero $f$\\
 \text{ \ } in the fraction field of $A$
such that $u' = fu$ and $v' = f^{-1}v$.
\end{propers}

\subsection*{An ad hoc tensor formalism}\text{} \\
 Let $A^{\otimes \mathbb{Z}}$
be a tensor product over $k$ of copies of $A$ indexed by
$\mathbb{Z}$ then $A^{\otimes \mathbb{z}}$ is a k-algebra
$A^{\otimes \mathbb{Z}}$ which is generated by symbols of
the form\\
\text{}\\
 \text{ \ } \ \text{ \ } \ \text{ \ } \ \text{ \ } \
\text{ \ } \ \text{ \ } \ \text{ \ } \ \text{ \ } \ $\bigotimes_{i
\in
\mathbb{Z}}(a_{i} \in A , a_{i}=1 \ for \ |i| \gg 0)$\\
\text{}\\
subject to the obvious relations of tensor algebra, under the
multiplication rule induced by componentwise multiplication in $A$\\
\text{} \\
\text{ \ } \ \text{ \ } \ \text{ \ } \ \text{ \ } \ \text{ \ } \
\text{ \ } \ \text{ \ } \ \text{ \ } \ $\bigotimes_{i \in \
\mathbb{Z}}a_{i} \cdot \bigotimes_{i \in \ \mathbb{Z}}b_{i}=
\bigotimes_{i \in \ \mathbb{Z}}a_{i}b_{i}$\\
\text{} \\
where for $i \in
\mathbb{Z}$, $a_{i},b_{i} \in A$\\
\text{} \\
Put\\
\text{ \ } \ \text{ \ } \ \text{ \ } \ $a^{(l)} \doteq \bigotimes_{i
\in \ \mathbb{Z}}\left\lbrace
           \begin{array}{c l}
              a & \text{if $i=l$ },\\
              1 & \text{otherwise}.
           \end{array}
         \right. $
for all $a \in A$ and $l \in \ \mathbb{Z}$.\\

For the universal property of $A^{\otimes \mathbb{Z}}$ see
\cite{And97} and \cite{Ger}.

%\begin{rem}{Universal construction of $A^{\otimes \mathbb{Z}}$\\}
%Let $\mathbf{A}$ be an $k-algebra$ equipped with a family $ \{ a
%\mapsto a^{(i)} : A  \rightarrow \mathbf{A} \}_{i= - \infty}^{
%\infty}$ of $k-algebra$ homomorphisms with the following universal
%property:\\
%For any $k-algebra$ $R$ and families $\{ \psi _{i} :A \rightarrow R
%\}_{i= - \infty}^{ \infty}$ of $k-algebra$ homomorphisms there exist
%a unique $k-algebra$ homomorphism $\Psi : \mathbf{A} \rightarrow R$
%such that $ \psi_{i}(a) = \Psi (a^{(i)})$ for all integers $i$ and
%$a \in A$.
%\\ Then $\mathbf{A}$ is unique up to isomorphism and $\mathbf{A}=A^{\otimes
%\mathbb{Z}}$.
%\end{rem}

\begin{defn}
Given a matrix $X$ with entries in $A$ and $l \in \mathbb{Z}$ define
$X^{(l)}$ to be a matrix with entries in $A^{\otimes \mathbb{Z}}$
such that
\[(X^{(l)})_{ij}=(X_{ij})^{(l)}.\]
\end{defn}

%\begin{rem}
%$ \{ a^{(l)} \}_{l \in \mathbb{Z}}$ generates \ $A^{\otimes
%\mathbb{Z}}$ as an $k-algebra$, due to the fact that  each element
%of
%\[\bigotimes_{i \in \mathbb{Z}}(a_{i} \in A , a_{i}=1 \ for \ |i|
%\gg 0) \] is an algebraic relation of a finite set of
%elements of $ \{ a^{(l)} \}_{l \in \mathbb{Z}}$ \\
%\end{rem}

For any subset $I \subset \ \mathbb{Z}$ let $A^{\otimes I}$ denote
the k-subalgebra of $A^{\otimes \mathbb{Z}}$ generated by all
elements of the form $a^{(l)}$ where $a \in A$ and $l \in I$,
moreover it is a universal commutative k-algebra with unit (i.e has
a universal property). If $I$ is a finite subset of $ \mathbb{Z}$,
then we can identify the k-algebra $A^{\otimes I}$ with the usual
tensor product over $k$ of copies
of $A$ indexed by $I$.\\

Note that by making the identifications $a^{(0)}=a$ for all $a \in
A$ we equip $A^{\bigotimes \mathbb{Z}}$ with the structure of an
A-algebra.

\begin{rem}
A well known fact is that the $k-algebra$ $A^{\bigotimes
\mathbb{I}}$ has no zero divisors.
\end{rem}

\begin{defn}
A \textbf{Jacobi matrix} $Z$ is an object admitting the following
properties
\begin{enumerate}
    \item \text{ $Z$ is an $n$ by $n$ matrix with entries in $A^{\bigotimes
    \{1,...,n+1\}}$.}
    \item $Z \neq 0$.
    \item $Z_{1 2} \in \text{$k$-span of} L \cdot L^{(1)}\cdot L^{(2)}\cdot
    (L^{(3)})^{\ 2} \cdots (L^{(n)})^{\ 2} \cdot L^{(n+1)}$.
    \item $Z_{1 1} = [1 \mapsto 2]_{\ast}Z_{1 2}$.
    \item $\pi_{\ast}Z_{i j}=Z_{\pi i , \pi j} \ \text{ for any
    bijective derangement $\pi$ supported in $\{ 1,...,n \}$}$.
    \item $\begin{vmatrix}
      Z_{1 1} & Z_{1 2} \\
      Z_{2 1} & Z_{2 2} \\
    \end{vmatrix} = 0$.
    \item $abel_{l=0}^{n+1}([0 \mapsto
-l]_{\ast}Z)=\bigtriangleup \cdot \bar{Z}$.\\
where $\bigtriangleup = \langle 2 | n+1\rangle_{\ast} \langle 0 |
1\rangle_{\ast}Z_{1 1} \cdot \langle 0 | 1\rangle_{\ast}Z_{1 1}
\cdot \langle 0 | 2\rangle_{\ast}Z_{2 2} \cdot
\prod_{l=3}^{n}{(\langle 0 | l\rangle_{\ast}Z_{l l})^{ \ 2}}$.\\
\end{enumerate}
\end{defn}

We call $\bigtriangleup$ the \textbf{\emph{discriminant}} of the
Jacobi matrix $Z$. If we need to draw attention to the basic data
we say that $Z$ is of \emph{type} $(k,n,A,L)$.\\

\section*{the abstract abel map}

\begin{defn}
Let $Y$ be an $n$ by $n$ matrix. The \textbf{abstract Abel map} is
defined to be $$ \textbf{Abel} : Y  \  \mapsto  \
abel_{l=0}^{n+1}Y^{(l)} $$

The \textbf{abstract Abel map} by definition sends each Segre matrix
$X$ to the $n$ by $n$ matrix $abel_{l=0}^{n+1}X^{(l)}$ with entries
in $A^{\bigotimes \{0,...,n+1 \}}$.

\end{defn}

\begin{prop}["The abstract Abel theorem"]\label{aat}
Let $X$ and $X'$ be Segre matrices with images $Z$ and $Z'$ under
the abstract Abel map, respectively. Then $X'$ is $k-equivalent$ to
$X$ if and only if $Z'$ is $k-proportional$ to $Z$.
\end{prop}

\begin{thm}{\text{}}\label{Jac7} \\
1. \ The set of $k-proportionality$ classes of $Jacobi-matrices$
forms a projective\\
 \text{ \ } \ variety.\\
2. \ The abstract Abel map takes values in the set of Jacobi
matrices.\\
3. \ The abstract Abel map puts the $k-equivalence$ classes of Segre
matrices into\\
\text{ \ } \ \ bijective correspondence with the $k-proportionality$
classes of $Jacobi-matrices$.
\end{thm}

\section*{G-forms}
\begin{defn}
Let $G$ be a divisor such that $$degG \equiv 0 mod 2 \ , \
\frac{1}{2}degG \geq 2g$$ and put $n = \frac{1}{2}degG - g + 1$.\\
A G-form $X$ is an object admitting the following properties:
\begin{enumerate}
    \item $X \text{ is an $n$ by $n$ matrix with entries in }
     \mathcal{L}(G)$.
    \item \text{All two by two minors of $X$ are zero}.
    \item There exists in $X$ some row and some column with
    $k-linearly$ independent entries.
\end{enumerate}
\end{defn}
A G-form is the same thing as a Segre matrix of type
$\begin{pmatrix}
  k , n , \bigoplus_{m=0}^{\infty}\mathcal{L}(mG) , \mathcal{L}(G) \\
\end{pmatrix}$\\
\begin{prop}{\textbf{Construction of a divisor represented by a given G-form.}}\label{gf4}
Let $G$ be a divisor such that $degG \equiv 0 mod 2 \ , \
\frac{1}{2}degG \geq 2g$ and $n = \frac{1}{2}degG - g + 1$. Let $X$
be a $G-form$. choose a column $u$ and a row $v$ of $X$ with
$k-linearly$ independent entries. Let $f$ be the entry common to $u$
and $v$; then $f$ does not vanish identically. Consider the
effective divisors
$$ D = G + min_{j}(div v_{j}), \ E = G + min_{i}(div u_{i}), F = G + divf.$$
Our claim is that $X$ represents the divisor $D$.
\end{prop}

\begin{prop}{\textbf{Representation of divisors of degree $\frac{1}{2}degG$ by
G-forms.}} Let $G$ be a divisor such that $degG \equiv 0 mod 2 \ , \
\frac{1}{2}degG \geq 2g$ and $n = \frac{1}{2}degG - g + 1$. Let $D$
be a divisor of degree $\frac{1}{2}degG$. Let $u$ be a column vector
with entries forming a $k-basis$ for $\mathcal{L}(D)$ and let $v$ be
a row vector with entries forming a $k-basis$ for $\mathcal{L}(G -
D)$. Put $X=uv$.\\
1. $X$ is a G-form.\\
2. The $k-equivalence$ class of $X$ depends only on $D$, not on the
choice of vectors\\
 \text{ \ } \ $u$ and $v$.
\end{prop}

\begin{defn}
In the above situation we say that the $G-form$ $X$
\emph{represents} the divisor $D$ of degree $\frac{1}{2}degG$.
\end{defn}

\begin{prop}\label{gf5}
Let $G$ be a divisor of even degree such that $\frac{1}{2}degG
\geq 2g$.\\
There exists a unique bijective correspondence $$\left\{%
\begin{array}{c}
    \text{k-equivalence classes of G-forms} \\
\end{array}
\right\} \leftrightarrow \left\{%
\begin{array}{c}
    \text{divisor classes of degree } \frac{1}{2}degG
 \\
\end{array}
\right \}.$$ with respect to which, for any $G-form$ $X$ and for any
divisor $D$ of degree $\frac{1}{2}degG$, the $k-equivalence$ class
of $X$ corresponds to the divisor class of $D$ if and only if $X$
represents $D$.
\end{prop}

\subsection*{Matrix representation of divisor class addition and
subtraction}

\begin{defn}{\textbf{Kronecker products.}}
Given a $p$ by $q$ matrix $A$ and an $r$ by $s$ matrix $B$ both with
entries in some ring $R$, the Kronecker product $A \circ B$ is
defined to be the $pr$ by $qs$ matrix with entries in $R$ admitting
a decomposition into $r$ by $s$ blocks of the form
$$ A \circ B = \begin{bmatrix}
       & \vdots &   \\
  \cdots & A_{i j} B & \cdots \\
       & \vdots &   \\
\end{bmatrix}.$$
\end{defn}

\begin{rem}
The \textbf{Kronecker product} of matrices is compatible with the
ordinary matrix multiplication in the sense that $$(A \circ B)(X
\circ Y)=( AX ) \circ ( BY )$$ whenever $AX$ and $BY$ are defined.
\end{rem}
\begin{prop}\label{mrep2}
Let divisors $G, G', D$ and $D'$ be given subject to the following
conditions:
$$degG = 2 \cdot degD,  degG' = 2 \cdot degD'$$
$$min \left(\begin{array}{c}
          \frac{1}{2} degG , \frac{1}{2} degG' \\
      \end{array} \right ) \geq 2g, \ max \left( \begin{array}{c}
          \frac{1}{2} degG , \frac{1}{2} degG' \\
      \end{array} \right ) \geq 2g + 1.$$
Put \text{ \ } \ \text{ \ } $n = \frac{1}{2} degG -g +1, \ n' =
\frac{1}{2} degG'
-g +1,$\\
and \text{ \ } \ \text{ \ } $n'' = \frac{1}{2} (degG +degG') -g +1 =
n + n' + g - 1.$\\ Fix a $G-form$ $X$ representing $D$ and a
$G'-form$ $X'$ representing $D'$. Let $P$ and $Q$ be any $nn'$ by
$nn'$ permutation matrices and consider the block decomposition
$$P (X \circ X') Q = \begin{bmatrix}
  a & b \\
  c & d \\
\end{bmatrix}$$
where the block $d$ is $n''$ by $n''$ and the other blocks are of
the appropriate sizes.\\
1. \ For some $P$ and $Q$  the corresponding block $d$ is
$k-general$.\\
2. \ For any $P$ and $Q$ such that $d$ is $k-general$, $d$ is a $(G
+ G')-form$ representing $D + D'$.
\end{prop}
\begin{lem}\label{mrep3}
Let E be a nonzero effective divisor. Let $\mathcal{R}_{E}$ be the
ring consisting of the meromorphic functions on $\mathrm{C}$ regular
in a neighborhood of the support of $E$ and let $\mathcal{I}_{E}
\subset \mathcal{R}_{E}$ be the ideal consisting of functions
vanishing to the order at least $E$(i.e. functions which have zeroes
at each point of the support of $E$, with multiplicities greater or
equal to the order of $E$ there). Then there exists a $k-linear$
functional
$$\sigma : \mathcal{R}_{E} \rightarrow k$$ factoring through the
quotient $\mathcal{R}_{E} / \mathcal{I}_{E}$ such that the induced
$k-bilinear$ map
$$
((a mod \mathcal{I}_{E} , b mod \mathcal{I}_{E} \mapsto \sigma
(ab)): \mathcal{R}_{E} / \mathcal{I}_{E} \times \mathcal{R}_{E} /
\mathcal{I}_{E} \rightarrow k
$$
is a perfect pairing of $(degE)$-dimensional vector spaces over $k$.
\end{lem}
\begin{lem}\label{mrep4}
Let $G$ and $E$ be divisors such that
$$ degG \equiv 0 mod 2, \ E > 0, \ \frac{1}{2}degG - degE > 2g - 2. $$
There exists a $k-linear$ functional $$\rho : \mathcal{L}(G)
\rightarrow k$$ factoring through the quotient
$\frac{\mathcal{L}(G)}{\mathcal{L}(G - E)}$ such that for all
divisors $D$ of degree $\frac{1}{2}degG$ the induced $k-bilinear$
map $$ \left( (a + \mathcal{L}(D - E) , b + \mathcal{L}(G - D - E))
\mapsto \rho(ab)\right) : \frac{\mathcal{L}(D)}{\mathcal{L}(D - E)}
\times \frac{\mathcal{L}(G - D)}{\mathcal{L}(G - D - E)} \rightarrow
k$$ is a perfect pairing of $(degE)$-dimensional vector spaces over
$k$.
\end{lem}
\text{ \ }\\
In the situation of Lemma \ref{mrep4} we call $\rho : \mathcal{L}(G)
\mapsto k$ an $E-compression functional$.
\begin{prop}\label{cmp}
Let $G$ and $E$ be divisors such that
$$ degG \equiv 0 mod 2, \ E > 0, \ \frac{1}{2}degG - degE \geq 2g$$
and put
$$n=\frac{1}{2}degG - degE - g + 1, \ n'=\frac{1}{2}degG - g + 1 = n + degE$$
Let $\rho : \mathcal{L}(G) \rightarrow k$ be an $E$-compression
functional. Let $D$ be a divisor of degree $\frac{1}{2}degG$ and let
$X$ be a $G-form$ representing $D$. Let $P$ and $Q$ be any $n'$ by
$n'$ permutation matrices and consider the block decomposition
$$P X Q = \begin{bmatrix}
  a & b \\
  c & d \\
\end{bmatrix}$$
where the block $a$ is $degE$ by $degE$, the block $d$ is $n$ by
$n$ and the other blocks are of the appropriate sizes.\\
1 . \ For some $P$ and $Q$ we have $| \rho a | \neq 0$.\\
2 . \ For any $P$ and $Q$ such that $| \rho a | \neq 0$ the matrix
$z$ defined by the rule
$$\begin{bmatrix}
  w & x \\
  y & z \\
\end{bmatrix} = \begin{bmatrix}
  1 & 0 \\
  -(\rho c)(\rho a)^{ \ -1} & 1 \\
\end{bmatrix}\begin{bmatrix}
  a & b \\
  c & d \\
\end{bmatrix}\begin{bmatrix}
  1 & -(\rho a)^{ \ -1}(\rho b) \\
  0 & 1 \\
\end{bmatrix}$$
is a $(G - 2E)-form$ representing $D - E$.
\end{prop}
\section*{Candidates for the Jacobian and for the Abel map}
\subsection{Candidate for the Jacobian}\text{ \ }\\
Fix an effective divisor $E$ of degree $ \geq 2g + 1$ and put
$$A=\bigoplus_{m=0}^{\infty} \mathcal{L}( 2mE ), \
n = \ell(E)=degE - g + 1, \ L = \mathcal{L}(2E).$$ The projective
algebraic variety $J$ of $k-proportionality$ classes of
$Jacobi-matrices$ of type $(k, n, A, L)$ is our candidate for the
Jacobian of $\mathcal{C}$.

\subsection*{Candidate for the Abel map}\text{ \ }\\
For each divisor $D$ of degree zero fix a $2E-form$
 \ $X_{D}$ representing $D + E$. Now a $2E-form$ is the same thing
 as a Segre matrix of type $(k, n, A, L)$. By Proposition \ref{gf5} it
 follows that the map $D \mapsto X_{D}$ puts the divisor classes
 of degree zero in bijective correspondence with the
 $k-equivalence$ classes of Segre matrices of type $(k, n, A, L)$.
 For each divisor $D$ of degree zero let $Z_{D}$ be the image of
 $X_{D}$ under the abstract Abel map. By Theorem \ref{Jac7} it follows
 that the map $ D \mapsto Z_{D} $ puts the classes of divisors of
 degree zero into bijective correspondence with the points of $J$.
 The bijective map from classes of divisors of degree zero to $J$
 induced by the map $D \mapsto Z_{D}$ is our candidate for the
 Abel map.

\medskip

\begin{lem}\label{cc2}
Fix an $E-compression$ functional $ \rho(4E) \rightarrow k$. Fix
Segre matrices $X$ and $X'$ of type $(k, n, A, L)$. Fix a divisor
$D$ (resp. $D'$) such that $X$ (resp. $X'$) is $k-equivalent$ to
$X_{D}$ (resp. $X_{D'}$). Let $P$ and $Q$ be any $n^{ \ 2}$ by $n^{
\ 2}$permutation matrices and consider the block decomposition
$$P (X \circ X')Q=\begin{bmatrix}
  \bullet & \bullet & \bullet \\
  \bullet & a & b \\
  \bullet & c & d \\
\end{bmatrix}$$
where the block $a$ is of size $degE$ by $degE$, the block $d$ is of
size $n$ by $n$ and the other blocks are of the appropriate sizes.
Here the bullets hold places for blocks the contents of which do not
concern us. Further, consider the block-decomposed matrix
$$\begin{bmatrix}
  w & x \\
  y & z \\
\end{bmatrix} =| \rho a | \cdot \begin{bmatrix}
  | \rho a | & 0 \\
  -(\rho c)(\rho a)^{\star} & | \rho a | \\
\end{bmatrix}\begin{bmatrix}
  a & b \\
  c & d \\
\end{bmatrix}\begin{bmatrix}
  | \rho a | & -(\rho a)^{\star}(\rho b) \\
  0 & | \rho a | \\
\end{bmatrix}$$
1. \ For some $P$ and $Q$ the corresponding block $z$ is
$k-general$.\\
2. \ For any $P$ and $Q$ such that the corresponding block $z$ is
$k-general$, $z$ is a Segre matrix of type $(k, n, A, L)$ and
moreover $z$ is $k-equivalent$ to $X_{D + D'}$
\end{lem}

\smallskip

\begin{proof}
Since $X$ (resp. $X'$) is $k-equivalent$ to $X_{D}$ (resp. $X_{D'}$)
they differ only by the $k-basis$ elements of $\mathcal{L}(E+D),
\mathcal{L}(2E - (E+D))$ (resp. $\mathcal{L}(E+D'), \ \mathcal{L}(2E
- (E+D'))$), thus we may consider $X_{D}$and $X_{D'}$.

We hereby list some facts to make the context of the proof more
vivid :\\
\textbf{A}. $X_{D}$ is a $2E-form$ representing $E + D$.\\
\textbf{B}. $X_{D}$ is of size $\ell(E + D) = deg(E + D) - g + 1 =
deg E+
   degD - g + 1 = degE - g + 1 = n$.\\
\textbf{C}. $X_{D} \circ X_{D'} = u_{D}v_{D} \circ u_{D'}v_{D'} =
(u_{D} \circ u_{D'})(v_{D} \circ v_{D'})$
   where $u_{D}$ (resp. $v_{D}$) is a column (resp. row) vector
   with entries forming a $k-basis$ for $\mathcal{L}(E+D)$ (resp.
   $\mathcal{L}(2E - (E+D))=\mathcal{L}(E-D)$) and so is $u_{D'}$ (resp. $v_{D'}$)
   replacing $D$ by $D'$.
   $u_{D} \circ u_{D'}$ is an $n^{ \ 2}$ column vector with
   entries $k$-spanning $\mathcal{L}((E + D)+(E + D'))=\mathcal{L}(2E + D +
   D')$, $v_{D} \circ v_{D'}$ is an $n^{ \ 2}$ row vector with
   entries $k$-spanning $\mathcal{L}((2E - (E+D))+(2E - (E+D)))=\mathcal{L}(4E - (2E + D +
   D'))=\mathcal{L}(2E - (D + D'))$. \\
\textbf{D}. By (\ref{mrep2}) we can choose permutation matrices
$P_{1}$ and
   $Q_{1}$ such that the above matrix multiplication gives us
   a $4E-form$ representing $(E + D)+(E + D')$ in the right-bottom  $\ell(2E + D +
   D')$ by $\ell(2E - (D + D'))$  block, where $\ell(2E + D +
   D') = \ell(2E - (D + D')) = deg2E -g +1 = degE + n $.\\
\textbf{E}. Using the $E-compression$ functional $ \rho(4E)
\rightarrow k$
   on that right-bottom block we decompose the $\mathcal{L}((E + D)+(E +
   D'))$ $k-basis$ to a $k-basis$ of $\frac{\mathcal{L}(2E + D + D')}{\mathcal{L}(E + D + D')}$
   and a $k-basis$ of $\mathcal{L}(E + D + D')$, and the $\mathcal{L}(2E - (D + D'))$ $k-basis$ to a $k-basis$ of $\frac{\mathcal{L}(2E - (D + D'))}{\mathcal{L}(E - (D + D'))}$
   and a $k-basis$ of $\mathcal{L}(E - (D + D'))$ .\\
\textbf{F}. Following Proposition \ref{cmp}, one more time we use a
permutation matrices $P_{2}$and $Q_{2}$ to bring the $k-basis$ of
$\frac{\mathcal{L}(2E + D + D')}{\mathcal{L}(E + D + D')}$
   to the first top $degE$ elements of the corresponding column
   vector
   (which is in fact the column vector formed by the last $degE + n$
   entries of $P_{1}(u_{D} \circ u_{D'})$) and to bring the $k-basis$ of $\frac{\mathcal{L}(2E - (D + D'))}{\mathcal{L}(E - (D + D'))}$
   to the first left $degE$ elements of the corresponding row
   vector
   (which is in fact the row vector formed by the last $degE + n$
   entries of $(v_{D} \circ v_{D'})Q_{1}$).
   Thus we have at the left bottom $n$ by $n$ block, $d$,of the
   referred  matrix a $2E-form$ representing $E + D + D'$\\
    In view of the above observations by taking $P=P_{2}P_{1}$ and

   $Q=Q_{1}Q_{2}$ as permutation matrices we have that
    $$\begin{bmatrix}
  a & b \\
  c & d \\
\end{bmatrix}$$ is a $4E-form$ representing $2E+D+D'$.\\
$| \rho a | \neq 0$. \\
$z$ is a 2E-form representing $E+D+D'$.\\
 A fortiori $z$ is
$k-general$.\\

   2.\ By applying $\rho$ to both sides of the above equation we get
that $$|\rho w|=|\rho a|^{ \ 3 }\rho a, \ \rho x = 0, \ \rho y =
0.$$ If $|\rho a|=0$ we immediately  have the right hand side of the
equation vanishing, thus we get that the matrix on the left is a
zero matrix particularly $z$ which contradicts its $k-generality$
so we have $|\rho a| \neq 0$ and hence $|\rho w| \neq 0$.\\
By hypothesis, since $X$ and $X'$ are Segre matrices we have a
factorization $X \circ X' = uv \circ u'v' = (u \circ u')(v \circ
v')$ .Since $X$ and $X'$ are $k-equivalent$ to $X_{D}$ and $X_{D'}$
this above factorization has the form
$$\begin{bmatrix}
  w & x \\
  y & z \\
\end{bmatrix} =\begin{bmatrix}
   p \\
   q \\
\end{bmatrix}\begin{bmatrix}
  r & s\\
\end{bmatrix},$$

where the entries of the column vector (resp. row vector) on the
right belong to $\mathcal{L}(2E + D + D')$ (resp. $\mathcal{L}(2E -
D - D')$),the blocks $p$ and $r$ are vectors of length $degE$ and
the blocks $q$ and $s$ are vectors of length $n$. By the definition
of $E$-compression functional (Lemma \ref{cmp}) it follows that the
entries of $p$ (resp. $r$) project to a $k-basis$ of the quotient
$\frac{\mathcal{L}(2E + D + D')}{\mathcal{L}(E + D + D')}$ (resp.
$\frac{\mathcal{L}(2E - (D + D'))}{\mathcal{L}(2E - (E + D + D'))}$)
since otherwise we would have had a zero row in $\rho (pr) = \rho a$
contradicting the fact that $|\rho a|$ does not vanish. Also by the
definition of the $E$-compression functional it follows that the
entries of $q$ (resp. $s$) belong to $\mathcal{L}(E + D + D')$
(resp. $\mathcal{L}(2E - (E + D + D'))$). Finally, since $z=qs$ is
$k-general$, the entries of $q$ (resp. $s$) must be $k-linearly$
independent, and hence the entries of $q$ (resp. $s$) must form a
$k-basis$ of $\mathcal{L}(E + D + D')$ (resp. $\mathcal{L}(2E - (E +
D + D'))$). Therefore the block $z$ is indeed a $2E-form$
representing $E + D + D'$ and hence $k-equivalent$ to $X_{D + D'}$.
\end{proof}

% ------------------------------------------------------------------------

\section*{example 1}

For some background on elliptic functions see Whittaker and Watson
\cite{Whit}, Chapter XX.\\
\smallskip

 We hereby introduce a part of the
method in the case of a nonsingular
projective plane cubic, i.e., the well known Weierstrass curve.\\
An affine part of this projective curve can be identified with the
space $\mathbb{C}/\Lambda$ via the map $$z \mapsto (1, \rho(z), \rho
'(z))$$ described in the background chapter of elliptic functions,
which in turn can be identified with a complex torus manifold. Since
any rational function $F \in k(\mathrm{C})$ can be written in the
form $f(1,x,y)$ where $x=\frac{w_{1}}{w_{0}}, y=\frac{w_{2}}{w_{0}}$
we have the map $f(1,x,y) \mapsto g(z)=f(1, \rho(z), \rho '(z))$
where $(1,x,y) \in \mathrm{C}$ inducing an isomorphism between the
rational functions field of the projective curve and the field of
meromorphic functions on $\mathbb{C}/\Lambda$ which
is just the field of elliptic functions. \\
The Weierstrass curve is the simplest nontrivial case in which the
Jacobian variety is nonzero group. In this case the Jacobian variety
is actually isomorphic to the curve itself, i.e, each divisor class
in $Pic^{0}(\mathrm{C})$ corresponds to a unique point on the curve,
namely it can be shown that fixing a point, say $b$ in $\mathrm{C}$
then the map $$\mathrm{C} \rightarrow Pic^{0}(\mathrm{C})$$ given by
$$t \mapsto [t-b]$$ defines the bijective correspondence. In view of
this observation we have that the
curve itself could be viewed as an additive group. \\
By the above notations, we hereby take $E=n[0]$ , $G=2E=2n[0]$, ,
 where $n \geq 2$ and $t \in \mathbb{C}/\Lambda$.
Since $degG=2n$, $degE=n$ and $g=1$ we have
$\frac{1}{2}degG > 2g$ and $degE=\frac{1}{2}degG$ like required.\\
Let $D$ be a divisor of degree zero representing a divisor class in
$Pic^{0}(\mathrm{C})$. Without loss of generality taking $b=[0]$ we
can take $D = [t]-[0]$ where $t \in \mathrm{C}$, i.e., $t \in
\mathbb{C}/\Lambda$. Hence we have $E + D = (n-1)[0] + [t]$, $G - (E
+ D) = E - D =(n+1)[0] - [t]$. We now have to find a
$\mathbb{C}-basis$ of $\mathcal{L}(E + D)$.Consider a rational
function $g \in k(\mathcal{C})$.

Under the above identification we have that $g \in
\mathcal{L}((n-1)[0] + [t])$ if and only if $\tilde{g} \in
\mathbb{C}/\Lambda$ is a meromorphic elliptic function such that an
the (not necessarily elliptic) function $f = \tilde{g} \cdot
\sigma(z - t)\sigma(z)^{n - 1}$ has no poles i.e., is entire complex
function. Hence we can identify $\mathcal{L}(E + D)$ with the space
of entire functions $f(z)$ such that the meromorphic function
$\frac{f(z)}{\sigma(z - t)\sigma(z)^{n - 1}}$ is $\Lambda
-periodic$. Hence in any given fundamental domain $f(z)$ must have
exactly $n$ zeros which sum up to $-t \in \mathbb{C}/\Lambda$ thus
$f(z)$ is of the form $f(z)=c \cdot \sigma(z-P_{1}) \cdots
\sigma(z-P_{n})$ where $P_{1}+\cdots+P_{n}=-t$ and $c \in
\mathbb{C}$. Now, the entries of the row vector
\footnotesize
$$\vec{\sigma}(z + \frac{t}{n})=
\begin{bmatrix} \sigma(z - \frac{t}{n})^{n}
  & \sigma(z - \frac{t}{n})^{n}\rho(z - \frac{t}{n})  & \sigma(z - \frac{t}{n})^{n}\rho '(z - \frac{t}{n})  & \cdots
    & \sigma(z - \frac{t}{n})^{n}\rho(z - \frac{t}{n})^{(n-2)} \\
    \end{bmatrix}$$\normalsize are $\mathbb{C}$ linearly independent entire functions
admitting these restrictions. Since we have by Riemann-Roch that the
space of such entire functions $f$ is an $n$-dimensional space, it
follows that these functions form a $\mathbb{C}$-basis the space of
the entire functions $f$. It follows that the entries of the row
vector $\frac{\vec{\sigma}(z + \frac{t}{n})}{\sigma(z -
t)\sigma(z)^{n - 1}}$ form a $\mathbb{C}$-basis for $g \in
\mathcal{L}((n-1)[0] + [t])$. In much of the same way we can
identify $\mathcal{L}(E - D)$ with the space of entire functions
$f(z)$ such that the meromorphic function $\frac{f(z)\sigma(z -
t)}{\sigma(z)^{n + 1}}$ is $\Lambda -periodic$. Thus $f(z)$ is of
the form $f(z)=c \cdot \sigma(z-Q_{1}) \cdots \sigma(z-Q_{n})$ where
$Q_{1}+\cdots+Q_{n}=t$ and $c \in \mathbb{C}$, and the entries of
the row vector $\frac{\sigma(z - t)\vec{\sigma}(z -
\frac{t}{n})}{\sigma(z)^{n + 1}}$ form a $\mathbb{C}$-basis for $g
\in \mathcal{L}((n+1)[0] - [t])$. Now, since $(\frac{\sigma(z -
t)\vec{\sigma}(z - \frac{t}{n})}{\sigma(z)^{n + 1}})^{T} \cdot
\frac{\vec{\sigma}(z + \frac{t}{n})}{\sigma(z - t)\sigma(z)^{n -
1}}=\frac{(\vec{\sigma}(z - \frac{t}{n}))^{T} \cdot \vec{\sigma}(z +
\frac{t}{n})}{\sigma(z)^{2n}}$, we have that $\frac{(\vec{\sigma}(z
- \frac{t}{n}))^{T} \cdot \vec{\sigma}(z +
\frac{t}{n})}{\sigma(z)^{2n}}$ is a $2E-form$ representing $D +
E$.\\
\text{}\\
We have\\
\text{}\\
 $\begin{pmatrix}
   abel_{l=0}^{n+1}(\vec{\sigma}(z - \frac{t}{n}))^{T} \cdot \vec{\sigma}(z +
\frac{t}{n}) \\
 \end{pmatrix}_{i j}=$\\
 \text{}\\
$\begin{vmatrix} & \vec{\sigma}(z_{0} + \frac{t}{n})
\\  & \vdots \\ & \vec{\sigma}(z_{i} +
\frac{t}{n}) \\ & \vdots \\\end{vmatrix}\begin{vmatrix}  \cdots
\vec{\sigma}(z_{i} - \frac{t}{n})  \cdots  \vec{\sigma}(z_{n+1} -
\frac{t}{n})\\ \end{vmatrix} \cdot \begin{vmatrix} & \vdots\\ &
\vec{\sigma}(z_{i} + \frac{t}{n}) \\ & \vdots
\\ & \vec{\sigma}(z_{n+1} +
\frac{t}{n})  \\\end{vmatrix}
\begin{vmatrix}
   \vec{\sigma}(z_{0} - \frac{t}{n})  \cdots  \vec{\sigma}(z_{j} - \frac{t}{n})  \cdots \\
\end{vmatrix}$
where
$$\begin{vmatrix} & \vec{\sigma}(z_{1})
\\  & \vdots \\ & \vec{\sigma}(z_{n}) \\\end{vmatrix}=\begin{vmatrix} \sigma(z_{1})^{n}
  & \sigma(z_{1})^{n}\rho(z_{1})  & \sigma(z_{1})^{n}\rho '(z_{1})  & \cdots
    & \sigma(z_{1})^{n}\rho(z_{1})^{(n-2)} \\ &  & \cdots  &   &
    \\ \sigma(z_{n})^{n}
    & \sigma(z_{n})^{n}\rho(z_{n})  & \sigma(z_{n})^{n}\rho '(z_{n})  & \cdots
    & \sigma(z_{n})^{n}\rho(z_{n})^{(n-2)} \\\end{vmatrix}=$$
    $$\prod_{i=0}^{n+1}\sigma(z_{i})^{n} \begin{vmatrix} 1
  & \rho(z_{1})  & \rho '(z_{1})  & \cdots
    & \rho(z_{1})^{(n-2)} \\ &  & \cdots  &   & \\ 1
    & \rho(z_{n})  & \rho '(z_{n})  & \cdots
    & \rho(z_{n})^{(n-2)} \\\end{vmatrix} =$$
$$\sum_{k=2}^{n}(-1)^{k-2}\cdot 0!1! \cdots (n)! \cdot
\sigma(\sum_{i=1}^{n}z_{i}) \cdot \prod_{1 \leq i \leq j \leq
n}\sigma(z_{i}-z_{j})$$ thus we have that
$$\begin{vmatrix} & \vdots\\ & \hat{z_{k}}
\\ & \hat{z_{m}} \\ & \vdots \\\end{vmatrix}=0!1! \cdots (n)! \cdot
\sigma(\sum_{i \neq k, m}z_{i}) \cdot \prod_{\begin{array}{c}
                           0 \leq i \leq j \leq n+1\\
                                 i \neq k, m    \\
                                             \end{array}
}\sigma(z_{i}-z_{j}).$$ Put $z_{i} \mapsto z_{i} + \frac{t}{n}$
then \\
$ \begin{vmatrix} \vdots \\ \vec{\sigma}(z_{k} + \frac{t}{n})
\\ \vdots \\ \vec{\sigma}(z_{m} +
\frac{t}{n}) \\  \vdots \\\end{vmatrix}= 0!1! \cdots (n)! \cdot
\sigma(\sum_{i \neq k, m}z_{i}+ \frac{t}{n}) \cdot
\prod_{\begin{array}{c}
                           0 \leq i \leq j \leq n+1\\
                                 i \neq k, m    \\
                                             \end{array}
}\sigma(z_{i}+ \frac{t}{n}-(z_{j}+ \frac{t}{n}))=$ \text{ \ } \
\text{ \ } \ \text{ \ } \  \text{ \ } \ \text{ \ } \ $=0!1! \cdots
(n)! \cdot \sigma(t + \sum_{i \neq k, m}z_{i} ) \cdot
\prod_{\begin{array}{c}
                           0 \leq i \leq j \leq n+1\\
                                 i \neq k, m    \\
                                             \end{array}
}\sigma(z_{i}-z_{j})$\\ symmetrically done with $z_{i} \mapsto z_{i}
- \frac{t}{n}$
\section*{example 2}
 Let us consider the following nonsingular curve of genus $g=3$ in
$\mathbb{P}^{2}$
$$C :  \ \omega_{2}^{4}=\omega_{1}(\omega_{1}-\omega_{0}))(\omega_{1}-2\omega_{0})(\omega_{1}-3\omega_{0})$$
where $\omega_{0}, \omega_{1}, \omega_{2}$ are homogeneous
coordinates of $\mathbb{P}^{2}$. \\
Consider the divisor $G=14P$, where $P=(1, 1, 0)\in \mathbb{P}^{2}$
then $degG=14=4g+2$. Consider the affine neighborhood $A_{0}=\{
\omega_{0} \neq 0 \}$. Let $x = \frac{\omega_{1}}{\omega_{0}}$ and
$y = \frac{\omega_{2}}{\omega_{0}}$ be the affine coordinates of
$A_{0}$, then the equation defining the curve $C$ in $A_{0}$ is
$$C : \  y^{4}=x(x-1)(x-2)(x-3)$$
and the affine representation of $P$ there is $(1,0)$.\\
In $\mathcal{O}_{C, P}$ we have $(x-1)=y_{4} \cdot
\frac{1}{x(x-2)(x-3)}$, since $x(x-2)(x-3)$ does not vanish at $P$
and so is an invertible element in the local ring, thus we have that
the (only) maximal ideal in  $\mathcal{O}_{C, P}$ is $m_{C, P}=y
\mathcal{O}_{C, P}$ and $y$ is the local parameter there. Another
way for obtaining this fact is by considering the gradient $\nabla=
(f'(x),4y^{3})$ where $f=x(x-1)(x-2)(x-3)$ at $(1,0)$ which is
$(f'(1),0)$, since $f$ is separable and $f(1)=0$ we have that $f'(1)
\neq 0$ then $y$ is a local parameter in a
neighborhood of $(1,0)$.\\
By the Riemann - Roch theorem we have $dimG = degG - g + 1 =
14-3+1=12$. Thus we have to find $12$ basis elements for
$\mathcal{L}(G)$.\\
$f_{1}(x,y)=\frac{1}{x-1}$, then in homogeneous coordinates we have
$\tilde{f_{1}}(\omega_{0}, \omega_{1}, \omega_{2})=(\omega_{0}
,\omega_{1}- \omega_{0})$. Now $\tilde{f}(\omega_{0}, \omega_{1},
\omega_{2})=(1,0)$ yields that $\omega_{0} = \omega_{1}$
intersecting with $C$ we have that $ \omega_{2} = 0$. Since the
multiplicity of $ \omega_{2}$ in the defining equation of the curve
is $4$ we obtain that $\tilde{f_{1}}$ has an order 4 pole at $P$ and
no other poles on the curve. $\tilde{f_{1}}(\omega_{0}, \omega_{1},
\omega_{2})=(0,1)$ yields that $\omega_{0}=0$ thus intersecting with
$C$ we have $\omega_{2}^{4}=\omega_{1}^{4}$ then $\tilde{f_{1}}$
gets four distinct zeros $Q_{1}=(0,1,1) ,
Q_{2}=(0,1,i) , Q_{3}=(0,1,-1)$ and  $Q_{4}=(0,1,-i)$.\\
Summing up our results we have
$$divf_{1} = Q_{1} + Q_{2} + Q_{3} + Q_{4} - 4P$$
thus $\tilde{f_{1}} \in \mathcal{L}(G)$.\\
Next we consider the function $f_{2}(x,y)=\frac{y}{x-1}$, another
representation for $f_{2}$ in the neighborhood of $P$ is
$f_{2}(x,y)=\frac{y \cdot x(x-2)(x-3) }{y^{4}}$. In homogenous
coordinates we have $\tilde{f_{2}}(\omega_{0}, \omega_{1},
\omega_{2})=(\omega_{1}(\omega_{1}-2\omega_{0})(\omega_{1}-3\omega_{0})
, \omega_{2}^{3})$.  Now $\tilde{f}(\omega_{0}, \omega_{1},
\omega_{2})=(1,0)$ yields that $ \omega_{2} = 0$ intersecting with
$C$ we get  $\omega_{1}=0$, $\omega_{1}=\omega_{0}$,
$\omega_{1}=2\omega_{0}$ and $\omega_{1}=3\omega_{0}$ thus we have
four points $(1,0,0), (1,1,0), (1,2,0)$ and $(1,3,0)$ respectively.
The only point which maps to $(0,1)$ is $(1,1,0)$ since all other
options make the first coordinate of the map $\tilde{f_{2}}$ vanish.
Since the multiplicity of $\omega_{2}$ in the second coordinate
defining $\tilde{f_{2}}$ is three we have that $\tilde{f_{2}}$ has a
single pole of order 3 at $P$. For finding the zeros of
$\tilde{f_{2}}$ let us return to the first representation
$f_{2}(x,y)=\frac{y}{x-1}$, thus we have $\tilde{f_{2}}(\omega_{0},
\omega_{1}, \omega_{2})=(\omega_{2} ,\omega_{1}- \omega_{0})$.
$\tilde{f_{1}}(\omega_{0}, \omega_{1}, \omega_{2})=(0,1)$ yields
that $\omega_{2}=0$ thus intersecting with $C$ we get four distinct
points $Q'_{1}=(1,0,0) , Q'_{2}=(1,1,0) , Q'_{3}=(1,2,0)$ and
$Q'_{4}=(1,3,0)$. As $Q'_{2}$ is not sent to $(0,1)$ we have $3$
distinct zeros :
$Q'_{1}=(1,0,0) , Q'_{3}=(1,2,0)$ and $Q'_{4}=(1,3,0)$.\\
Summing up our results we have
$$divf_{2} = Q'_{1} + Q'_{3} + Q'_{4} - 3P$$
thus $\tilde{f_{2}} \in \mathcal{L}(G)$.\\

We now build the remaining functions as follows :\\
$f_{3}(x,y)=\frac{1}{(x-1)^{2}}$ thus $f_{3}=(f_{1})^{2}$ hence
$divf_{3} = 2 \cdot divf_{1} = 2Q_{1}
+2Q_{2} + 2Q_{3}+ 2Q_{4} - 8P.$ \\
$f_{4}(x,y)=\frac{y}{(x-1)^{2}}$ thus $f_{4}=f_{1} \cdot f_{2}$
hence $divf_{4} = divf_{1} + divf_{2} = Q_{1} + Q_{2} + Q_{3}
+ Q_{4} + $\\
\text{ \ } \ $ Q'_{1} + Q'_{3} + Q'_{4} - 7P.$\\
$f_{5}(x,y)=\frac{y^{2}}{(x-1)^{2}}$ thus $f_{5}=(f_{2})^{2}$ hence
$divf_{5} = 2 \cdot divf_{2} = 2Q'_{1} +
2Q'_{3} + 2Q'_{4} - 6P$\\
$f_{6}(x,y)=\frac{1}{(x-1)^{3}}$ thus $f_{3}=(f_{1})^{3}$ hence
$divf_{3} = 3 \cdot divf_{1} = 3Q_{1} + 3Q_{2} + 3Q_{3} + 3Q_{4} -
12P.$\\
$f_{7}(x,y)=\frac{y}{(x-1)^{3}}$ thus $f_{7}=(f_{1})^{2} \cdot
f_{2}$ hence $divf_{7} = 2 \cdot divf_{1} + divf_{2} = 2Q_{1}
+ 2Q_{2} + 2Q_{3} + $ \text{ \ } \ $2Q_{4} + Q'_{1} + Q'_{3} + Q'_{4} - 11P.$ \\
$f_{8}(x,y)=\frac{y^{2}}{(x-1)^{3}}$ thus $f_{8}=f_{1} \cdot
(f_{2})^{2}$ hence $divf_{8} = divf_{1} + 2 \cdot divf_{2} = Q_{1} +
Q_{2} + Q_{3} +$\\ \text{ \ } \ $ Q_{4} + 2Q'_{1} + 2Q'_{3} + 2Q'_{4} - 10P.$\\
$f_{9}(x,y)=\frac{y^{3}}{(x-1)^{3}}$ thus $f_{9}=(f_{2})^{3}$ hence
$divf_{9} = 3 \cdot divf_{2} = 3Q'_{1} + 3Q'_{3} + 3Q'_{4} -
9P$\\
$f_{10}(x,y)=\frac{y^{2}}{(x-1)^{4}}$ thus $f_{10}=(f_{1})^{2} \cdot
(f_{2})^{2}$ hence $divf_{10} = 2 \cdot divf_{1} + 2 \cdot divf_{2}
=$ \\ \text{ \ } \  $2Q_{1} + 2Q_{2} + 2Q_{3} + 2Q_{4} + 2Q'_{1} +
2Q'_{3} + 2Q'_{4} -
14P.$\\
$f_{11}(x,y)=\frac{y^{3}}{(x-1)^{4}}$ thus $f_{11}= f_{1} \cdot
(f_{2})^{3}$ hence $divf_{11} = divf_{1} + 3 \cdot divf_{2} = Q_{1}
+ Q_{2} + Q_{3} +$\\ \text{ \ } \ $ Q_{4} + 3Q'_{1} + 3Q'_{3} + 3Q'_{4} - 13P.$\\

Finally we take $f_{0}$ to be constant function.\\
We have 12 functions $f_{k}$ \ $k=0, ... ,11$ where $f_{k} \in
\mathcal{L}(G)$. Since for $i \neq j$ $f_{i}$ and $f_{j}$ have
distinct number of poles we have that  $f_{k}$ \ $k=0, ... ,11$ are
linearly  independent over $\mathbb{C}$ and thus form our desired
basis for $\mathcal{L}(G)$.\\
Now we construct a $2E$-form ,where $degE=2g+1=7$. We take $E=7P$
thus by the above computations we have the basis
$$\{f_{0}, \cdots , f_{11} \}$$
for $2E=14P$. Consider the $2E$-form
$$X_{D}=\begin{pmatrix}
  f_{1} & f_{5} & f_{2} & f_{4} & f_{0} \\
  f_{7} & f_{11} & f_{8} & f_{10} & f_{4} \\
  f_{3} & f_{8} & f_{4} & f_{7} & f_{1} \\
  f_{8} & f_{2}^{4} & f_{9} & f_{11} & f_{5} \\
  f_{4} & f_{9} & f_{5} & f_{8} & f_{2} \\
\end{pmatrix}$$

by \ref{gf4} (with the current notations) we have
$$ (E + D) = G + min_{j}(div v_{j}).$$
$min_{j}(div v_{j})=min \{Q_{1} + Q_{2} + Q_{3} + Q_{4} + Q'_{1} +
Q'_{3} + Q'_{4} - 7P ,3Q'_{1} + 3Q'_{3} + 3Q'_{4} - 9P ,\\ 2Q'_{1} +
2Q'_{3} + 2Q'_{4} - 6P, Q_{1} + Q_{2} + Q_{3} + Q_{4} + 2Q'_{1} +
2Q'_{3} + 2Q'_{4} - 10P, Q'_{1} + Q'_{3} + Q'_{4} - 3P \}\\ = Q'_{1}
+ Q'_{3} + Q'_{4} - 10P$. \\So, \ $E + D = Q'_{1} + Q'_{3} + Q'_{4}
- 10P + 2E$.\\ Hence $$D = Q'_{1} + Q'_{3} + Q'_{4} -10P + E
\Rightarrow D = Q'_{1} + Q'_{3} + Q'_{4} - 3P$$ and $degD = 3 - 3 =
0$, thus this $2E$-form represents the divisor $D$.

Note that in fact this $2E$-form represents the divisor $E + D =
Q'_{1} + Q'_{3} + Q'_{4} + 4P$. \\
 Let us take
$g_{1}(x,y)=\frac{x}{x-1}$, then in homogeneous coordinates we have
$\tilde{g_{1}}(\omega_{0}, \omega_{1}, \omega_{2})=(\omega_{1}
,\omega_{1}- \omega_{0})$. Now $\tilde{g}(\omega_{0}, \omega_{1},
\omega_{2})=(1,0)$ yields that $\omega_{0} = \omega_{1}$
intersecting with $C$ we have that $ \omega_{2} = 0$. Since the
multiplicity of $ \omega_{2}$ in the defining equation of the curve
is $4$ we obtain that $\tilde{g_{1}}$ has an order 4 pole at $P$ and
no other poles on the curve. $\tilde{g_{1}}(\omega_{0}, \omega_{1},
\omega_{2})=(0,1)$ yields that $\omega_{1}=0$ thus intersecting with
$C$ we have $\omega_{2}^{4}=0$ then $\tilde{9_{1}}$ gets four
zero of order 4 at  $Q_{1}=(1,0,0)$.\\
Summing up our results we have
$$divg_{1} = 4Q'_{1} - 4P$$
thus $\tilde{g_{1}} \in \mathcal{L}(G)$.\\
Now, we define :\\
$g_{2}=f_{2}$.\\
$g_{3}=(g_{1})^{2}$  hence $divg_{3} = 2 \cdot divg_{1} = 8Q'_{1} -
8P.$\\
$g_{4}=g_{1} \cdot g_{2}$ hence $divg_{4} = divg_{1} + divg_{2} =
5Q'_{1} + Q'_{3} + Q'_{4} - 7P.$\\
$g_{5}=(g_{2})^{2}$ hence $divg_{5} = 2 \cdot divg_{2} = 2Q'_{1} +
2Q'_{3} + 2Q'_{4} - 6P$\\
$g_{6}=(g_{1})^{3}$ hence $divg_{6} = 3 \cdot divg_{1} = 12Q'_{1} -
12P.$\\
$g_{7}=(g_{1})^{2} \cdot g_{2}$ hence $divg_{7} = 2 \cdot divg_{1} +
divg_{2} = 9Q'_{1} + Q'_{3} + Q'_{4} - 11P.$\\
$g_{8}=g_{1} \cdot (g_{2})^{2}$ hence $divg_{8} = divg_{1} + 2 \cdot
divg_{2} = 6Q'_{1} + 2Q'_{3} + 2Q'_{4} - 10P.$\\
$g_{9}=(g_{2})^{3}$ hence $divg_{9} = 3 \cdot divg_{2} = 3Q'_{1} +
3Q'_{3} + 3Q'_{4} - 9P$ \\
$g_{10}=(g_{1})^{2} \cdot (g_{2})^{2}$ hence $divg_{10} = 2 \cdot
divg_{1} + 2 \cdot divg_{2} = 10Q'_{1} + 2Q'_{3} + 2Q'_{4} - 14P.$\\
$g_{11}= g_{1} \cdot (g_{2})^{3}$ hence $divg_{11} = divg_{1} + 3
\cdot divg_{2} = 7Q'_{1} + 3Q'_{3} + 3Q'_{4} - 13P.$
finally we take $g_{0}$ to be constant function.\\
\text{}\\
For the same considerations as for $f_{k}$, \ $k=0, ... ,11$, we
have that $g_{k}$, \ $k=0, ... ,11$ form a basis for
$\mathcal{L}(G)$.\\
\text{}\\
Consider the $2E$-form
$$X_{D'}=\begin{pmatrix}
  g_{2} & g_{5} & g_{1} & g_{4} & g_{0} \\
  g_{8} & g_{11} & g_{7} & g_{10} & g_{4} \\
  g_{5} & g_{9} & g_{4} & g_{8} & g_{2} \\
  g_{9} & g_{2}^{4} & g_{8} & g_{11} & g_{5} \\
  g_{4} & g_{8} & g_{3} & g_{7} & g_{1} \\
\end{pmatrix}$$
\text{}\\
By \ref{gf4} (with the current notations) we have
$$ (E + D') = G + min_{j}(div v_{j}).$$
$min_{j}(div v_{j})=min \{4Q'_{1}- 4P, 9Q'_{1} + Q'_{3} + Q'_{4} -
11P , 8Q'_{1} - 8P, 6Q'_{1} + 2Q'_{3} + 2Q'_{4} - 10P,\\ \text{ \ }
\ \text{ \ } \ \text{ \ } \ \text{ \ } \  5Q'_{1} + Q'_{3} + Q'_{4}
- 7P
\} = 4Q'_{1} -11P$. \\
So, \ $E + D' = 4Q'_{1} - 11P + 2E$.\\
Hence $$D' = 4Q'_{1} -11P + E \Rightarrow D = 4Q'_{1} - 4P$$ and
$degD' = 4 - 4 = 0$, thus this $2E$-form represents the divisor
$D'$.\\

We now perform the group operation. First note that\\
\small
$X_{D}=u_{D} \cdot v_{D}=\begin{pmatrix}
  f_{2}^{-1} \\
  f_{1} \\
  f_{1}f_{2}^{-1} \\
  f_{2} \\
  f_{0} \\
\end{pmatrix} \cdot \begin{pmatrix}
  f_{1}f_{2} \\
  (f_{2})^{3} \\
  (f_{2})^{2} \\
  f_{1}(f_{2})^{2} \\
  f_{2} \\
\end{pmatrix}^{T}$
, $X_{D'}=u_{D'} \cdot v_{D'} = \begin{pmatrix}
  g_{1}^{-1} \\
  g_{2} \\
  g_{2}g_{1}^{-1} \\
  (g_{2})^{2}g_{1}^{-1} \\
  g_{0} \\
\end{pmatrix} \cdot \begin{pmatrix}
  g_{1}g_{2} \\
  g_{1}(g_{2})^{2} \\
  (g_{1})^{2} \\
  g_{2}(g_{1})^{2} \\
  g_{1}\\
\end{pmatrix}^{T}$\\
\normalsize
 $X_{D} \circ X_{D'} =( u_{D} \cdot v_{D} ) \circ ( u_{D'} \cdot v_{D'} ) = (u_{D} \circ u_{D'}) \cdot (v_{D}
\circ v_{D'})=$ \\
\  \text{ \ } \ \text{ \ } \  \text{ \ } \ \text{ \ }\  \text{ \ } \
\text{ \ } \ \text{ \ } \ \text{ \ } \ \text{ \ } $\begin{pmatrix}
  f_{2}^{-1}g_{1}^{-1} \\
  1 \\
  g_{1}^{-1} \\
  g_{1}^{-1}g_{2} \\
  f_{2}^{-1}g_{0} \\
  f_{1}g_{1}^{-1} \\
  f_{1}g_{2} \\
  f_{1}g_{1}^{-1}g_{2} \\
  f_{1}g_{1}^{-1}(g_{2})^{2} \\
  f_{1}g_{0} \\
  f_{1}f_{2}^{-1}g_{1}^{-1} \\
  f_{1} \\
  f_{1}g_{1}^{-1} \\
  f_{1}g_{1}^{-1}g_{2} \\
  f_{1}f_{2}^{-1}g_{0} \\
  f_{2}g_{1}^{-1} \\
  g_{2}^{2} \\
  g_{1}^{-1}g_{2}^{2} \\
  g_{1}^{-1}(g_{2})^{3} \\
  f_{2}g_{0} \\
  f_{0}g_{1}^{-1} \\
  f_{0}g_{2} \\
  f_{0}g_{1}^{-1}g_{2} \\
  f_{0}g_{1}^{-1}(g_{2})^{2} \\
  f_{0}g_{0} \\
\end{pmatrix}^{T}  \cdot \begin{pmatrix}
  f_{1}g_{1}(g_{2})^{2} \\ f_{1}g_{1}(g_{2})^{3} \\ f_{1}f_{2}(g_{1})^{2} \\ f_{1}(g_{2})^{2}(g_{1})^{2} \\  f_{1}f_{2}g_{1}  \\
g_{1}(g_{2})^{4} \\ g_{1}(g_{2})^{5} \\ (f_{2})^{3}(g_{1})^{2} \\ (g_{1})^{2}(g_{2})^{4} \\  (f_{2})^{3}g_{1}  \\
g_{1}(g_{2})^{3} \\ g_{1}(g_{2})^{4} \\ (f_{2})^{2}(g_{1})^{2} \\ (g_{1})^{2}(g_{2})^{3} \\  (f_{2})^{2}g_{1}  \\
f_{1}g_{1}(g_{2})^{3} \\ f_{1}g_{1}(g_{2})^{4} \\ f_{1}(f_{2})^{2}(g_{1})^{2} \\ f_{1}(g_{1})^{2}(g_{2})^{3} \\  f_{1}(f_{2})^{2}g_{1}  \\
g_{1}(g_{2})^{2} \\ g_{1}(g_{2})^{3} \\ f_{2}(g_{1})^{2} \\ (g_{1})^{2}(g_{2})^{2} \\  f_{2}g_{1}  \\
   \end{pmatrix}$.\\
\text{}\\
\textbf{Note}: In the above expressions, whenever there where
instances of $f_{2}$ and $g_{2}$ in the same element we have written
them both as $g_{2}$ for simplicity.\\
\text{}\\

We will now find the elements in $u_{D} \circ u_{D'}$ which form a
basis for $\mathcal{L}(D + E + D'+ E) = \mathcal{L}(5Q'_{1} + Q'_{2}
+ Q'_{3} + 7P)$, and the elements in $v_{D} \circ v_{D'}$ which form
a basis for $\mathcal{L}(2E-(D + E) + 2E-(D'+ E)) = \mathcal{L}(2E -
D - D') = \mathcal{L}(21P - 5Q'_{1} - Q'_{2} - Q'_{3})$\\

First notice that $f_{1} = g_{1} - 1$ and $f_{2} = g_{2}$. Now,
applying these equalities and taking $1$ for $g_{0}$ and $f_{0}$ we
rewrite $u_{D} \circ u_{D'}$ as indicated in figure \textbf{ (a) }.
It is now very easy to see the linear dependence between the
elements. We have the following set of $12$ linear independent
elements: \\
$\{ 1 , g_{1}-1 , g_{1}^{-1} , g_{1}^{-1}g_{2} ,
g_{1}^{-1}(g_{2})^{2} , g_{1}^{-1}(g_{2})^{3} ,
g_{1}^{-1}g_{2}^{-1}, g_{1}g_{2}^{-1} - g_{2}^{-1} , g_{1}g_{2} -
g_{2} , g_{2} , g_{2}^{-1} , (g_{2})^{2} \}.$ \\
\newpage
All these functions belong to $\mathcal{L}(D + E + D'+ E) =
\mathcal{L}(5Q'_{1} + Q'_{3} + Q'_{4} + 7P)$. As $deg(D + E + D'+
E)=14$ and $g=3$ we have by Riemann-Roch $dim(D + E + D'+ E)= 14 - 3
+ 1 = 12$, thus indeed these functions form a basis for
$\mathcal{L}(D + E + D'+ E)$. By multiplying $(u_{D} \circ u_{D'})$
by a permutation matrix from its right, we can put these basis
elements in the bottom $12$ entries of the column vector.\\
Using the same procedure described above we get the vector in figure \textbf{ (b) }\\

\small
 $(u_{D} \circ u_{D'})=
\begin{pmatrix}
  g_{1}^{-1}g_{2}^{-1} \\
  1 \\
  g_{1}^{-1} \\
  g_{1}^{-1}g_{2} \\
  g_{2}^{-1} \\
  1 - g_{1}^{-1} \\
  g_{1}g_{2} - g_{2}\\
  g_{2} - g_{1}^{-1}g_{2}\\
  (g_{2})^{2} - g_{1}^{-1}(g_{2})^{2}\\
  g_{1} - 1\\
  g_{2}^{-1} - g_{1}^{-1}g_{2}^{-1}\\
  1 - g_{1}\\
  1 - g_{1}^{-1} \\
  g_{2} - g_{1}^{-1}g_{2} \\
  g_{1}g_{2}^{-1} - g_{2}^{-1}\\
  g_{2}g_{1}^{-1} \\
  g_{2}^{2} \\
  g_{1}^{-1}g_{2}^{2} \\
  g_{1}^{-1}(g_{2})^{3} \\
  g_{2} \\
  g_{1}^{-1} \\
  g_{2} \\
  g_{1}^{-1}g_{2} \\
  g_{1}^{-1}(g_{2})^{2} \\
  1 \\
\end{pmatrix}$ \ \textbf{ (a) } \text{ \ } \ \text{ \ } \  $(v_{D} \circ v_{D'})=
\begin{pmatrix}
  g_{1}^{2}g_{2}^{2} - g_{1}g_{2}^{2} \\
  g_{1}^{2}g_{2}^{3} - g_{1}g_{2}^{3} \\
  g_{1}^{3}g_{2} - g_{1}^{2}g_{2} \\
  g_{1}^{3}g_{2}^{2} - g_{1}^{2}g_{2}^{2} \\
  g_{1}^{2}g_{2} - g_{1}g_{2} \\
  g_{1}g_{2}^{4} \\
  g_{1}g_{2}^{5} \\
  g_{1}^{2}g_{2}^{3} \\
  g_{1}^{2}g_{2}^{4} \\
  g_{1}g_{2}^{3}\\
  g_{1}g_{2}^{3}\\
  g_{1}g_{2}^{4}\\
  g_{1}^{2}g_{2}^{2}\\
  g_{1}^{2}g_{2}^{3} \\
  g_{1}g_{2}^{2} \\
  g_{1}^{2}g_{2}^{3} - g_{1}g_{2}^{3}\\
  g_{1}^{2}g_{2}^{4} - g_{1}g_{2}^{4} \\
  g_{1}^{3}g_{2}^{2} - g_{1}^{2}g_{2}^{2} \\
  g_{1}^{3}g_{2}^{3} - g_{1}^{2}g_{2}^{3} \\
  g_{1}^{2}g_{2}^{2} - g_{1}g_{2}^{2} \\
  g_{1}g_{2}^{2} \\
  g_{1}g_{2}^{3} \\
  g_{1}^{2}g_{2} \\
  g_{1}^{2}g_{2}^{2} \\
  g_{1}g_{2} \\
\end{pmatrix}^{T}$ \text{ \ } \ \textbf{ (b) }\\ \normalsize
\\
\\
We have the following set of $12$ linear independent elements:\\
\\
\small$\{ g_{1}g_{2} , g_{1}g_{2}^{2} , g_{1}^{2}g_{2} ,
g_{1}g_{2}^{3} , g_{1}^{2}g_{2}^{2} , g_{1}^{3}g_{2} -
g_{1}^{2}g_{2} , g_{1}g_{2}^{4} , g_{1}^{2}g_{2}^{3} ,
g_{1}^{3}g_{2}^{2} - g_{1}^{2}g_{2}^{2} , g_{1}g_{2}^{5} ,
g_{1}^{2}g_{2}^{4} ,
g_{1}^{3}g_{2}^{3} - g_{1}^{2}g_{2}^{3} \}.$\\
\\
\normalsize
All these functions belong to $\mathcal{L}(2E-(D + E) +
2E-(D'+ E)) = \mathcal{L}(2E - D - D') = \mathcal{L}(21P - 5Q'_{1} -
Q'_{3} - Q'_{4})$. As $deg(2E - D - D')=14$ and $g=3$ we have by
Riemann-Roch $dim(2E - D - D')= 14 - 3 + 1 = 12$ thus indeed these
functions form a basis for $\mathcal{L}(2E - D - D')$. By
multiplying $(v_{D} \circ v_{D'})$ by a permutation matrix we can
put these basis
elements in the right $12$ entries of the row vector.\\

By multiplication by elementary matrices (with nonvanishing
determinant) we can assume that the last 12 bottom entries of
$(u_{D} \circ u_{D'})$ and the last 12 right entries of $(v_{D}
\circ v_{D'})$ have respectively the forms

$$\begin{pmatrix}
  1 \\
  g_{1} \\
  g_{2} \\
  g_{1}g_{2} \\
  (g_{2})^{2} \\
  g_{1}^{-1} \\
  g_{1}^{-1}g_{2} \\
  g_{1}^{-1}(g_{2})^{2} \\
  g_{1}^{-1}(g_{2})^{3} \\
  g_{1}^{-1}g_{2}^{-1} \\
  g_{1}g_{2}^{-1} \\
  g_{2}^{-1} \\
\end{pmatrix}=
\begin{pmatrix}
  1 \\
  \frac{\omega_{1}}{\omega_{1}- \omega_{0}} \\
  \frac{\omega_{2}}{\omega_{1}- \omega_{0}} \\
  \frac{\omega_{1}\omega_{2}}{(\omega_{1}- \omega_{0})^{2}} \\
  \frac{\omega_{2}^{2}}{(\omega_{1}- \omega_{0})^{2}} \\
  \frac{\omega_{1}- \omega_{0}}{\omega_{1}} \\
  \frac{\omega_{2}}{\omega_{1}} \\
  \frac{(\omega_{2})^{2}}{\omega_{1}(\omega_{1}- \omega_{0})} \\
  \frac{(\omega_{2})^{3}}{\omega_{1}(\omega_{1}- \omega_{0})^{2}} \\
  \frac{(\omega_{1}- \omega_{0})^{2}}{\omega_{1}\omega_{2}}  \\
  \frac{\omega_{1}}{\omega_{2}} \\
  \frac{\omega_{1}- \omega_{0}}{\omega_{2}} \\
\end{pmatrix} \  \text{ \ } \text{ \ }  \text{and} \text{ \ } \text{ \ }
\begin{pmatrix}
  g_{1}g_{2} \\
  g_{1}g_{2}^{2} \\
  g_{1}^{2}g_{2} \\
  g_{1}g_{2}^{3} \\
  g_{1}^{2}g_{2}^{2} \\
  g_{1}^{3}g_{2} \\
  g_{1}g_{2}^{4} \\
  g_{1}^{2}g_{2}^{3} \\
  g_{1}^{3}g_{2}^{2} \\
  g_{1}g_{2}^{5}\\
  g_{1}^{2}g_{2}^{4}\\
  g_{1}^{3}g_{2}^{3}\\
\end{pmatrix}^{T}=
\begin{pmatrix}
  \frac{\omega_{1}\omega_{2}}{(\omega_{1}- \omega_{0})^{2}} \\
  \frac{\omega_{1}\omega_{2}^{2}}{(\omega_{1}- \omega_{0})^{3}} \\
  \frac{\omega_{1}^{2}\omega_{2}}{(\omega_{1}- \omega_{0})^{3}} \\
  \frac{\omega_{1}\omega_{2}^{3}}{(\omega_{1}- \omega_{0})^{4}} \\
  \frac{\omega_{1}^{2}\omega_{2}^{2}}{(\omega_{1}- \omega_{0})^{4}} \\
  \frac{\omega_{1}^{3}\omega_{2}}{(\omega_{1}- \omega_{0})^{4}} \\
  \frac{\omega_{1}\omega_{2}^{4}}{(\omega_{1}- \omega_{0})^{5}} \\
  \frac{\omega_{1}^{2}\omega_{2}^{3}}{(\omega_{1}- \omega_{0})^{5}} \\
  \frac{\omega_{1}^{3}\omega_{2}^{2}}{(\omega_{1}- \omega_{0})^{5}} \\
  \frac{\omega_{1}\omega_{2}^{5}}{(\omega_{1}- \omega_{0})^{6}} \\
  \frac{\omega_{1}^{2}\omega_{2}^{4}}{(\omega_{1}- \omega_{0})^{6}} \\
  \frac{\omega_{1}^{3}\omega_{2}^{3}}{(\omega_{1}- \omega_{0})^{6}} \\
\end{pmatrix}^{T}$$

Multiplying these column and row vectors we get $d$:

$$\begin{pmatrix}
 g_{1}g_{2} & g_{1}g_{2}^{2} & g_{1}^{2}g_{2} & g_{1}g_{2}^{3} & g_{1}^{2}g_{2}^{2} & g_{1}^{3}g_{2} & g_{1}g_{2}^{4} & g_{1}^{2}g_{2}^{3} & g_{1}^{3}g_{2}^{2} & g_{1}g_{2}^{5} \ \ g_{1}^{2}g_{2}^{4} \ \ g_{1}^{3}g_{2}^{3} \\
  g_{1}^{2}g_{2} & g_{1}^{2}g_{2}^{2} & g_{1}^{3}g_{2} & g_{1}^{2}g_{2}^{3} & g_{1}^{3}g_{2}^{2} & g_{1}^{4}g_{2}^{2} & g_{1}^{2}g_{2}^{4} & g_{1}^{3}g_{2}^{3} & \mathbf{g_{1}^{4}g_{2}^{2}} & g_{1}^{2}g_{2}^{5} \ \ g_{1}^{3}g_{2}^{4} \ \ g_{1}^{4}g_{2}^{3} \\
  g_{1}g_{2}^{2} & g_{1}g_{2}^{3} & g_{1}^{2}g_{2}^{2} & g_{1}g_{2}^{4} & g_{1}^{2}g_{2}^{3} & g_{1}^{3}g_{2}^{2} & g_{1}g_{2}^{5} & g_{1}^{2}g_{2}^{4} & g_{1}^{3}g_{2}^{3} & \mathbf{g_{1}g_{2}^{6}} \ \ g_{1}^{2}g_{2}^{5} \ \ g_{1}^{3}g_{2}^{4} \\
  g_{1}^{2}g_{2}^{2} & g_{1}^{2}g_{2}^{3} & g_{1}^{3}g_{2}^{2} & g_{1}^{2}g_{2}^{4} & g_{1}^{3}g_{2}^{3} & \mathbf{g_{1}^{4}g_{2}^{2}} & g_{1}^{2}g_{2}^{5} & g_{1}^{3}g_{2}^{4} & g_{1}^{4}g_{2}^{3} & g_{1}^{2}g_{2}^{6} \ \ g_{1}^{3}g_{2}^{5} \ \ g_{1}^{4}g_{2}^{4} \\
  g_{1}g_{2}^{3} & g_{1}g_{2}^{4} & g_{1}^{2}g_{2}^{3} & g_{1}^{2}g_{2}^{5} & g_{1}^{2}g_{2}^{4} & g_{1}^{3}g_{2}^{3} & \mathbf{g_{1}g_{2}^{6}} & g_{1}^{2}g_{2}^{5} & g_{1}^{3}g_{2}^{4} & g_{1}g_{2}^{7} \ \ g_{1}^{2}g_{2}^{6} \ \ g_{1}^{3}g_{2}^{5} \\
  g_{2} & g_{2}^{2} & g_{1}g_{2} & g_{2}^{3} & g_{1}g_{2}^{2} & g_{1}^{2}g_{2} & g_{2}^{4} & g_{1}g_{2}^{3} & g_{1}^{2}g_{2}^{2} & g_{2}^{5} \ \ g_{1}g_{2}^{4} \ \ g_{1}^{2}g_{2}^{3} \\
  g_{2}^{2} & g_{2}^{3} & g_{1}g_{2}^{2} & g_{2}^{4} & g_{1}g_{2}^{3} & g_{1}^{2}g_{2}^{2} & g_{2}^{5} & g_{1}g_{2}^{4} & g_{1}^{2}g_{2}^{3} & g_{2}^{6} \ \ g_{1}g_{2}^{5} \ \ g_{1}^{2}g_{2}^{4} \\
  g_{2}^{3} & g_{2}^{4} & g_{1}g_{2}^{3} & g_{2}^{5} & g_{1}g_{2}^{4} & g_{1}^{2}g_{2}^{3} & g_{2}^{6} & g_{1}g_{2}^{5} & g_{1}^{2}g_{2}^{4} & g_{2}^{7} \ \ \mathbf{g_{1}g_{2}^{6}} \ \ g_{1}^{2}g_{2}^{5} \\
  g_{2}^{4} & g_{2}^{5} & g_{1}g_{2}^{4} & g_{2}^{6} & g_{1}g_{2}^{5} & g_{1}^{2}g_{2}^{4} & g_{2}^{7} & \mathbf{g_{1}g_{2}^{6}} & g_{1}^{2}g_{2}^{5} & g_{2}^{8} \ \ g_{1}g_{2}^{7} \ \ g_{1}^{2}g_{2}^{6} \\
  1 & g_{2} & g_{1} & g_{2}^{2} & g_{1}g_{2} & g_{1}^{2} & g_{2}^{3} & g_{1}g_{2}^{2} & g_{1}^{2}g_{2} & g_{2}^{4} \ \ g_{1}g_{2}^{3} \ \ g_{1}^{2}g_{2}^{2} \\
  g_{1}^{2} & g_{1}^{2}g_{2} & g_{1}^{3} & g_{1}^{2}g_{2}^{2} & g_{1}^{3}g_{2} & g_{1}^{4} & g_{1}^{2}g_{2}^{3} & g_{1}^{3}g_{2}^{2} & g_{1}^{4}g_{2} & g_{1}^{2}g_{2}^{4} \ \ g_{1}^{3}g_{2}^{3} \ \ \mathbf{g_{1}^{4}g_{2}^{2}} \\
  g_{1} & g_{1}g_{2} & g_{1}^{2} & g_{1}g_{2}^{2} & g_{1}^{2}g_{2} & g_{1}^{3} & g_{1}g_{2}^{3} & g_{1}^{2}g_{2}^{2} & g_{1}^{3}g_{2} & g_{1}g_{2}^{4} \ \ g_{1}^{2}g_{2}^{3} \ \ g_{1}^{3}g_{2}^{2} \\
\end{pmatrix}$$\\
\\
\small
$g_{1}g_{2}=\frac{\omega_{1}\omega_{2}}{(\omega_{1}-\omega_{0})^{2}}
\ \
g_{1}g_{2}^{2}=\frac{\omega_{1}\omega_{2}^{2}}{(\omega_{1}-\omega_{0})^{3}}
\ \
g_{1}^{2}g_{2}=\frac{\omega_{1}^{2}\omega_{2}}{(\omega_{1}-\omega_{0})^{3}}
\ \
g_{1}g_{2}^{3}=\frac{\omega_{1}\omega_{2}^{3}}{(\omega_{1}-\omega_{0})^{4}}
\ \
g_{1}^{3}g_{2}=\frac{\omega_{1}^{3}\omega_{2}}{(\omega_{1}-\omega_{0})^{4}}$\\
$g_{1}^{2}g_{2}^{2}=\frac{\omega_{1}^{2}\omega_{2}^{2}}{(\omega_{1}-\omega_{0})^{4}}
\ \
g_{1}g_{2}^{4}=\frac{\omega_{1}\omega_{2}^{4}}{(\omega_{1}-\omega_{0})^{5}}
\ \
g_{1}^{2}g_{2}^{3}=\frac{\omega_{1}^{2}\omega_{2}^{3}}{(\omega_{1}-\omega_{0})^{5}}
\ \
g_{1}^{3}g_{2}^{2}=\frac{\omega_{1}^{3}\omega_{2}^{2}}{(\omega_{1}-\omega_{0})^{5}}
\ \
g_{1}^{4}g_{2}=\frac{\omega_{1}^{4}\omega_{2}}{(\omega_{1}-\omega_{0})^{5}}$\\
$g_{1}g_{2}^{5}=\frac{\omega_{1}\omega_{2}^{5}}{(\omega_{1}-\omega_{0})^{6}}
\ \
g_{1}^{2}g_{2}^{4}=\frac{\omega_{1}^{2}\omega_{2}^{4}}{(\omega_{1}-\omega_{0})^{6}}
\ \
g_{1}^{3}g_{2}^{3}=\frac{\omega_{1}^{3}\omega_{2}^{3}}{(\omega_{1}-\omega_{0})^{6}}
\ \
g_{1}^{4}g_{2}^{2}=\frac{\omega_{1}^{4}\omega_{2}^{2}}{(\omega_{1}-\omega_{0})^{6}}
\ \
g_{1}^{5}g_{2}=\frac{\omega_{1}^{5}\omega_{2}}{(\omega_{1}-\omega_{0})^{6}}$\\
$g_{1}g_{2}^{6}=\frac{\omega_{1}\omega_{2}^{6}}{(\omega_{1}-\omega_{0})^{7}}
\ \
g_{1}^{2}g_{2}^{5}=\frac{\omega_{1}^{2}\omega_{2}^{5}}{(\omega_{1}-\omega_{0})^{7}}
\ \
g_{1}^{3}g_{2}^{4}=\frac{\omega_{1}^{3}\omega_{2}^{4}}{(\omega_{1}-\omega_{0})^{7}}
\ \
g_{1}^{4}g_{2}^{3}=\frac{\omega_{1}^{4}\omega_{2}^{3}}{(\omega_{1}-\omega_{0})^{7}}
\ \
g_{1}^{5}g_{2}^{2}=\frac{\omega_{1}^{5}\omega_{2}^{2}}{(\omega_{1}-\omega_{0})^{7}}$\\
$g_{1}^{6}g_{2}=\frac{\omega_{1}^{6}\omega_{2}}{(\omega_{1}-\omega_{0})^{7}}
\ \
g_{1}g_{2}^{7}=\frac{\omega_{1}\omega_{2}^{7}}{(\omega_{1}-\omega_{0})^{8}}
\ \
g_{1}^{2}g_{2}^{6}=\frac{\omega_{1}^{2}\omega_{2}^{6}}{(\omega_{1}-\omega_{0})^{8}}
\ \
g_{1}^{3}g_{2}^{5}=\frac{\omega_{1}^{3}\omega_{2}^{5}}{(\omega_{1}-\omega_{0})^{8}}
\ \
g_{1}^{4}g_{2}^{4}=\frac{\omega_{1}^{4}\omega_{2}^{4}}{(\omega_{1}-\omega_{0})^{8}}$\\
$g_{1}^{5}g_{2}^{3}=\frac{\omega_{1}^{5}\omega_{2}^{3}}{(\omega_{1}-\omega_{0})^{8}}
\ \
g_{1}^{6}g_{2}^{2}=\frac{\omega_{1}^{6}\omega_{2}^{2}}{(\omega_{1}-\omega_{0})^{8}}
\ \
g_{1}^{7}g_{2}=\frac{\omega_{1}^{7}\omega_{2}}{(\omega_{1}-\omega_{0})^{8}}$\\
\normalsize
\\
\newpage
 Now the bottom $12$ entries of $u_{D} \circ u_{D'}$ and the
right $12$ entries of the row vector $(v_{D} \circ v_{D'})$ form a
$12 \times 12$ block $d$ on the right bottom corner of the $25
\times 25$ matrix $X_{D} \circ X_{D'}$(this is just the block $d$ in
Proposition \ref{mrep2}). The block $d$ is a $4E$-form
representing the divisor $D + D'+ 2E$.\\

\bigskip\goodbreak

We will now build the $E$-compression functional $ \rho :
\mathcal{L}(4E)
\rightarrow k$.\\
First, by Lemma \ref{mrep4}, we have to choose a meromorphic
function $f_{4E}$ on $\mathcal{C}$ by the following rule:\\

\  \text{ \ } \ \text{ \ } \  \text{ \ } \ \text{ \ } \ \text{ \ } \
\text{ \ } \ \text{ \ } \  \text{ \ } \ \text{ \ } \ \text{ \ } \
\text{ \ } \ $ord_{x}f_{4E} = ord_{x}4E$ \\
for all points $x \in suppE \cap supp4E = suppE$,\\

\  \text{ \ } \ \text{ \ } \  \text{ \ } \ \text{ \ } \ \text{ \ } \
\text{ \ } \ \text{ \ } \  \text{ \ } \ \text{ \ } \ \text{ \ } \
\text{ \ } \ $ord_{x}f_{4E} = 0$\\
for all points $x \in suppE \setminus
supp4E = \emptyset $.\\

We thus choose the function $f_{4E}=(g_{1}^{-1})^{7}$ where
$g_{1}=\frac{\omega_{1}}{\omega_{1}- \omega_{0}}$. As we have shown
above $divg_{1} = 4Q'_{1} - 4P$ where $Q'_{1}=(1,0,0)$ and
$P=(1,1,0)$. Thus $divf_{4E} = -7 \cdot (4Q'_{1} - 4P) = 28P -
28Q'_{1}$ and $f_{4E}$ fulfils the required property.\\
Now according to Lemma \ref{mrep3} we have to find a meromorphic
differential $\omega$ on $\mathcal{C}$ such that\\
\  \text{ \ } \ \text{ \ } \  \text{ \ } \ \text{ \ } \ \text{ \ } \
\text{ \ } \ \text{ \ } \  \text{ \ } \ \text{ \ } \ \text{ \ } \
\text{ \ } \ $ord_{x}\omega + ord_{x}E=0$\\
for all points $x \in supp E$, i.e, \\
\  \text{ \ } \ \text{ \ } \  \text{ \ } \ \text{ \ } \ \text{ \ } \
\text{ \ } \ \text{ \ } \  \text{ \ } \ \text{ \ } \ \text{ \ } \
\text{ \ } \ $ord_{P}\omega = -ord_{P}E = -7.$\\

Considering the affine neighborhood $A_{0}=\{ \omega_{0} \neq 0 \}$
we define $\omega = \frac{xy dy}{(x-1)^{2}}$, i.e, $\omega =
g_{1}(x,y) \cdot g_{2}(x,y)dy$ where
$x=\frac{\omega_{1}}{\omega_{0}}$ and
$y=\frac{\omega_{2}}{\omega_{0}}$; as we have shown above $y$ is a
local parameter and $ord_{P}\omega=-7$. Now, by Lemma \ref{mrep4} we
have \\ \small $\rho(4E)(h)=\sigma(f_{4E} h)=Res_{P}( h
\frac{(x-1)^{7}}{x^{7}}\frac{xydy}{(x-1)^{2}})=Res_{P}(h
\frac{(x-1)^{5}ydy}{x^{6}} )=Res_{P}(h
\frac{y^{21}dy}{(x-2)^{5}(x-3)^{5}x^{11}} ).$ \\\normalsize As we
calculate the residue in an analytic neighborhood of $P \in A_{0}$
we can take $h$ to be the affine representation of a function
$\tilde{h} \in
\mathcal{L}(4E)$.\\
In our case one can compute the cases where the residue is not zero
just by looking at the divisors of the functions in the $12 \times
12$ matrix. We are looking for functions with pole divisor equal to
$22P$. Only such functions will have a simple pole at $P$
 after applying the $E$-compression functional $ \rho(4E)$ (since $-22 + 28 - 7 = -1$, where $28 = ord_{P}f_{4E}$ and
$-7=ord_{P}\omega$). There are exactly two functions which admit
this condition,
$g_{1}^{4}g_{2}^{2}=\frac{\omega_{1}^{4}\omega_{2}^{2}}{(\omega_{1}-\omega_{0})^{6}}$
and
$g_{1}g_{2}^{6}=\frac{\omega_{1}\omega_{2}^{6}}{(\omega_{1}-\omega_{0})^{7}}$.
 where $div(g_{1}^{4}g_{2}^{2})= 18Q'_{1} + 2Q'_{3} + 2Q'_{4} -
 22P$ and $div(g_{1}g_{2}^{6})= 10Q'_{1} + 6Q'_{3} + 6Q'_{4} -
 22P$. We have highlighted these entries of the matrix above.
We have found $degE = deg7P = 7$ such entries. Since these entries
occur at distinct rows and columns, using row and column permutation
matrices (the matrices $P_{2}$ and $Q_{2}$ in Lemma \ref{cc2}) we
can bring these entries to form a block in the upper left corner of
the matrix in which they will be placed on the diagonal. Thus using
the notation of Lemma \ref{cc2} we get the matrix
$$P (X \circ X')Q=\begin{bmatrix}
  \bullet & \bullet & \bullet \\
  \bullet & a & b \\
  \bullet & c & d \\
\end{bmatrix}$$
where the block $a$ is of size $degE$ by $degE$, i.e, $7 \times 7$,
the block $d$ is of size $n$ by $n$ where $n=\frac{1}{2}degG - degE
- g + 1 = \frac{28}{2} - 7 - 3 + 1 = 5$ , where the entries of $d$
become zeroes after applying the $E$-compression functional $
\rho(4E)$ and the other blocks
are of the appropriate sizes.\\
Notice, the block $d$ here is \textbf{not} the block considered
previously in the construction.\\
Looking at the  entries of the row and column vectors used to create
the last matrix we see that row and column vector creating
\textbf{the block} $d$ are the column vector entries $$\{ g_{1}g_{2}
, g_{1}g_{2}^{2} , g_{1}^{2}g_{2} , g_{1}g_{2}^{3} ,
g_{1}^{2}g_{2}^{2}  \}$$ and the row vector entries $$\{ 1 ,
g_{1}^{-1} , g_{1}^{-1}g_{2} , g_{1}^{-1}g_{2}^{-1}, g_{2}^{-1}
\}.$$Each of these sets contains linearly independent elements as
observed above. The first set forms a basis for $\mathcal{L}((2E + D
+ D') - E) = \mathcal{L}(E + (D + D')) = \mathcal{L}(5Q'_{1} +
Q'_{3} + Q'_{4})$ and the second set forms a basis for
$\mathcal{L}(4E - (2E + D + D') - E) = \mathcal{L}(2E - (E + (D +
D'))).$\\
Thus using the notation of Lemma \ref{cc2} the block $z$ is a
$2E$-form representing $E + (D + D')$, i.e, by agreement (see
Candidate for the Abel map section) $z = X_{D + D'}$ represents $(D
+ D')$ as desired.

This is the place to note that, looking at the  elements of the row
(resp. column) vector used to create the last matrix,there are
exactly two options for an element's column (resp. row) in the
matrix to behave under $ \rho$. The first is that all the column's
(resp. row's) entries vanish and then it is guaranteed that the
element belongs to $\mathcal{L}(4E - (2E + D + D') - E)$
 (resp. $\mathcal{L}((2E + D + D') - E)$), the second is that there is
an entry that does not vanish under $ \rho$ in which case it is
guaranteed that the element belongs to $\frac{\mathcal{L}(4E - (2E +
D + D'))}{\mathcal{L}(4E - (2E + D + D') - E}$
(resp.$\frac{\mathcal{L}(2E + D + D')}{\mathcal{L}((2E + D + D') -
E)}$ ). Thus, as our non vanishing entries appearing in distinct
rows and columns of the matrix, they point out exactly the
separation of interest both in the row and the column vector.

\newpage

\section*{Summary and conclusions}

Summing up our results we can point out a number of issues.\\
\ \\
$\bullet$ First and most important, the algorithm works. Example
number $2$ performs the algorithm on a genus $3$ curve which is even
not hyperelliptic. As far as we know, explicit algorithms for this
case are still poorly studied.
 \\

\medskip\goodbreak
 \ \\
$\bullet$ For cryptographic implementations there is much interest
in further generalizing this algorithm to work on finite fields.
Analyzing the algorithm for future generalization, we consider two
distinct cases: the field $k$ is infinite and not algebraically
closed, the field $k$ is finite. We have spotted two issues that
must be taken in consideration. The first involves the explicit
representation of the Jacobian by a variety of Jacobi matrices. The
Abel map (the map providing the representation) properties are
strongly based on the fact that the $k-algebra$ $A^{\bigotimes
\mathbb{I}}$ (where $I$ is a finite set contained in the set of
integers) is a domain. This result is true if $k$ is an
algebraically closed field, otherwise this claim is generally not
true. The second issue involves the implementation of the group
operation (performed on the $G$-forms). More specifically, the
problem is the existence of a particular function $f_G$ \cite{And02}
which gives rise to the compression functional. Under the relevant
assumptions, the existence of
such a function is guaranteed only over infinite fields.\\

We conclude that in the case where $k$ is an infinite field (not
algebraically closed), the explicit representation obtained in this
work does not apply directly but the group law still works (on
G-forms). In the case where $k$ is a finite field, both issues are
to be investigated.

% ------------------------------------------------------------------------
%Included for Gather Purpose only:
%input "Xbib.bib"
\bibliographystyle{amsplain}
\bibliography{xbib}
\end{document}